\documentclass[12pt,draft]{amsart}

\usepackage[cp1251]{inputenc}
\usepackage[T2A]{fontenc}
\usepackage[english]{babel}
\usepackage{amsmath,amsfonts,amssymb}
\usepackage{geometry}
\newtheorem{theorem}{Theorem}
\newtheorem{corollary}{Corollary}

\newtheorem{proposition}{Proposition}

\theoremstyle{remark}

\theoremstyle{definition}

\author{Vladimir A. Mikhailets, Alexandr A. Murach}
\title[Elliptic systems in a refined scale]{Elliptic systems of
pseudodifferential equations in a refined scale on a closed
manifold}

\address{Institute of Mathematics NAS of Ukraine \\
    Tereshchenkivska str., 3 \\
    Kyiv\\
    Ukraine\\
    01601}

\email{mikhailets@imath.kiev.ua, murach@imath.kiev.ua}

\subjclass[2000]{Primary 35J45, Secondary 46E35}

\date{14/11/2007}

\keywords{Elliptic system, pseudodifferential operator, regularly
varying function, scale of spaces, the H\"ormander spaces, the
Fredgolm property}

\begin{document}

\maketitle

\begin{abstract}
We study a system of pseudodifferential equations that is elliptic
in the sense of Petrovskii on a closed compact smooth manifold. We
prove that the operator generated by the system is Fredholm one on
a refined two-sided scale of the functional Hilbert spaces.
Elements of this scale are the special isotropic spaces of
H\"{o}rmander--Volevich--Paneah.
\end{abstract}

\section{Introduction} In this article we consider
an elliptic in Petrovskii's sense system of linear
pseudodifferential equations on a closed smooth manifold. It is
well known (see e.g. \cite{Her87, Agr90}\,) that the operator $A$
corresponding to this system is bounded and Fredholm in
appropriate pairs of the Sobolev spaces. We investigate this
operator in the Hilbert scale of the special isotropic
H\"ormander--Volevich-- Paneah spaces [3--6]
$$
H^{s,\varphi}:=H_{2}^{\langle\cdot\rangle^{s}\,\varphi(\langle\cdot\rangle)},
\quad\langle\xi\rangle:=\bigl(1+|\xi|^{2}\bigr)^{1/2}. \eqno(1)
$$
Here, $s\in\mathbb{R}$ and $\varphi$ is a functional parameter
slowly varying at $+\infty$ in Karamata's sense. In particular,
every standard function
$$
\varphi(t)=(\log t)^{r_{1}}(\log\log
t)^{r_{2}}\ldots(\log\ldots\log
t)^{r_{n}},\quad\{r_{1},r_{2},\ldots,r_{n}\}\subset\mathbb{R},\;n\in\mathbb{N},
$$
is admissible. This scale was introduced and investigated by the
authors in \cite{MM2a, MM2b}. It contains Sobolev's scale
$\{H^{s}\}\equiv\{H^{s,1}\}$ and is attached to it by the number
parameter $s$ and being considerably finer.

Spaces of form (1) arise naturally in different spectral problems:
convergence of spectral expansions of self-adjoint elliptic
operators almost everywhere, in the norm of the spaces $L_{p}$
with $p>2$ or $C$ (see survey \cite{AIN76}); spectral asymptotics
of general self-adjoint elliptic operators in a bounded domain,
the Weyl formula, a sharp estimate of the remainder in it (see
\cite{Mikh82, Mikh89}) and others. They may be expected to be
useful in other "fine"\, questions.  Due to their interpolation
properties, the spaces $H^{s,\varphi}$ occupy a special position
among the spaces of a generalized smoothness which are actively
investigated and used today (see survey \cite{KalLiz87}, recent
articles \cite{HarMou04, FarLeo06} and the bibliography given
there).

The main result of this article is the theorem on boundedness and
the Fredholm property of the operator $A$ in scale (1). The
refined local smoothness of a solution of the elliptic system is
obtained as a significant application. Also some auxiliary results
which may be of interest by themselves are given. The case of
scalar differential operators was investigated earlier in [8,
15--18].

\section{The statement of the problem} Let $\Gamma$ be a closed
(compact and without a boundary) infinitely smooth manifold of
dimension $n\geq1$. We suppose that a certain $C^{\infty}$-density
$dx$ is defined on $\Gamma$. By $\mathcal{D}'(\Gamma)$ we denote
the linear topological space of all distributions on $\Gamma$,
that is $\mathcal{D}'(\Gamma)$ is a space antidual to the space
$C^{\infty}(\Gamma)$ with respect to the extension of the scalar
product in $L_{2}(\Gamma,dx)$ by continuity. This extension is
denoted by $(f,w)_{\Gamma}$ for $f\in\mathcal{D}'(\Gamma)$, $w\in
C^{\infty}(\Gamma)$.

We consider a system of linear equations
$$
\sum_{k=1}^{p}\:A_{j,k}\,u_{k}=f_{j}\;\;\mbox{on}\;\;\Gamma,\;\;
j=1,\ldots,p. \eqno(2)
$$
Here, $p\in\mathbb{N}$ and $A_{j,k}$, $j,k=1,\ldots,p$, are scalar
classical (polyhomogeneous) pseudodifferential operators of
arbitrary real orders defined on the manifold $\Gamma$ (see e.g.
\cite[Sec. 2.1]{Agr90}). A complete symbol of the
pseudodifferential operator $A_{j,k}$ is an infinitely smooth
complex-valued function on the cotangent bundle $T^{\ast}\Gamma$.
A principal symbol of $A_{j,k}$ which is positively homogeneous of
order $\mathrm{ord}\,A_{j,k}$ in every section
$T^{\ast}_{x}\Gamma\setminus\{0\}$, $x\in\Gamma$ and, moreover, is
not identically equal to zero, is also defined. We consider
equations (2) in the sense of the distribution theory, so
$u_{k},\,f_{j}\in\mathcal{D}'(\Gamma)$. We put for every index
$k=1,\ldots,p$
$$
m_{k}:=\max\left\{\mathrm{ord}\,A_{1,k},\ldots,\mathrm{ord}\,A_{p,k}\right\}.
$$

Let us assume system (2) to be \textit{elliptic in Petrovskii's
sense }, that is
$$
\det\left(\,a^{(0)}_{j,k}(x,\xi)\,\right)_{j,k=1}^{p}\neq0
\quad\mbox{for each}\;\;x\in\Gamma,\;\xi\in
T^{\ast}_{x}\Gamma\setminus\{0\}.
$$
Here $a_{j,k}^{(0)}(x,\xi)$ is the principal symbol of the
pseudodifferential operator $A_{j,k}$ in the case
$\mathrm{ord}\,A_{j,k}=m_{k}$, or $a_{j,k}^{(0)}(x,\xi)\equiv0$ in
the case $\mathrm{ord}\,A_{j,k}<m_{k}$.

Let us rewrite system (2) in the matrix form: $Au=f$ on $\Gamma$,
where $A:=\bigl(\,A_{j,k}\,\bigr)$\; is a square matrix of order
$p$, and $u=\mathrm{col}\,(u_{1},\ldots,u_{p})$,
$f=\mathrm{col}\,(f_{1},\ldots,f_{p})$ are functional columns. The
mapping $u\mapsto Au$ is a linear continuous operator in the space
$\bigl(\mathcal{D}'(\Gamma)\bigr)^{p}$.

\section{A refined scale of spaces} We denote by $\mathcal{M}$ the
set of all Borel measurable functions
$\varphi:[1,+\infty)\rightarrow(0,+\infty)$ such that the
functions $\varphi$ and $1/\varphi$ are bounded on every closed
interval $[1,b]$, where $1<b<+\infty$, and the function $\varphi$
is slowly varying at $+\infty$ in Karamata's sense, that is
$$
\lim_{t\rightarrow\,+\infty}{\varphi(\lambda\,t)}/{\varphi(t)}=1
\quad\mbox{for each}\quad\lambda>0.
$$

Let $s\in\mathbb{R}$, $\varphi\in\mathcal{M}$. We denote by
$H^{s,\varphi}(\mathbb{R}^{n})$  the set of all tempered
distributions $u$ such that the Fourier transform $\widehat{u}$ of
the distribution $u$ is a function locally Lebesgue integrable in
$\mathbb{R}^{n}$ which satisfies the condition
$$
\int\langle\xi\rangle^{2s}\,\varphi^{2}(\langle\xi\rangle)
\:|\widehat{u}(\xi)|^{2}\,d\xi<\infty.
$$
Here the integral is evaluated over $\mathbb{R}^{n}$, and
$\langle\xi\rangle:=(1+\xi_{1}^{2}+\ldots+\xi_{n}^{2})^{1/2}$. In
the space $\mathrm{H}^{s,\varphi}(\mathbb{R}^{n})$ we define the
inner product
$$
\bigl(u,v\bigr)_{\mathrm{H}^{s,\varphi}(\mathbb{R}^{n})}:=
\int\langle\xi\rangle^{2s}\varphi^{2}(\langle\xi\rangle)
\,\widehat{u}(\xi)\,\overline{\widehat{v}(\xi)}\,d\xi.
$$

The space $H^{s,\varphi}(\mathbb{R}^{n})$ is a special isotropic
Hilbert case of the spaces introduced by L. H\"ormander \cite[Sec.
2.2]{Her65}, \cite[Sec. 10.1]{Her86} and L. R. Volevich, B. P.
Paneah \cite[Sec. 2]{VoPa65}, \cite[Sec. 1.4.2]{Pa00}. In the
simplest case where $\varphi(\cdot)\equiv1$ the space
$H^{s,\varphi}(\mathbb{R}^{n})$ coincides with the Sobolev space.
It follows from the inclusions
$$
\bigcup_{\varepsilon>0}H^{s+\varepsilon}(\mathbb{R}^{n})=:H^{s+}(\mathbb{R}^{n})
\subset H^{s,\varphi}(\mathbb{R}^{n})\subset
H^{s-}(\mathbb{R}^{n}):=\bigcap_{\varepsilon>0}H^{s-\varepsilon}(\mathbb{R}^{n})
$$
that in the collection of spaces
$$
\{H^{s,\varphi}(\mathbb{R}^{n}):s\in\mathbb{R},\varphi\in\mathcal{M}\,\}
\eqno(3)
$$
the functional parameter $\varphi$ defines an additional
(subpower) smoothness with respect to the basic (power)
$s$-smoothness. Otherwise speaking, $\varphi$ \textit{refines} the
power smoothness.

The \textit{refined scale} over the manifold $\Gamma$ is
constructed from scale (3) in the usual way. Let us take a finite
atlas from the $C^{\infty}$-structure on $\Gamma$ consisting of
the local charts $\alpha_{j}:\mathbb{R}^{n}\leftrightarrow U_{j}$,
$j=1,\ldots,r$. Here the open sets $U_{j}$ form the finite
covering of the manifold $\Gamma$. Let functions $\chi_{j}\in
C^{\infty}(\Gamma)$, $j=1,\ldots,r$, form a partition of unity on
$\Gamma$ satisfying the condition $\mathrm{supp}\,\chi_{j}\subset
U_{j}$.

Let us set
$$
H^{s,\varphi}(\Gamma):=\left\{h\in\mathcal{D}'(\Gamma):\;
(\chi_{j}h)\circ\alpha_{j}\in
H^{s,\varphi}(\mathbb{R}^{n})\;\;\mbox{for
every}\;\;j=1,\ldots,r\right\}.
$$
Here $(\chi_{j}h)\circ\alpha_{j}$ is the representation of the
distribution $\chi_{j}h$ in the local chart $\alpha_{j}$. The
inner product in the space $H^{s,\varphi}(\Gamma)$ is defined by
the formula
$$
(f,g)_{H^{s,\varphi}(\Gamma)}:=\sum_{j=1}^{r}\,((\chi_{j}f)\circ\alpha_{j},
(\chi_{j}\,g)\circ\alpha_{j})_{H^{s,\varphi}(\mathbb{R}^{n})}
$$
and induces the norm
$\|h\|_{s,\varphi}:=\bigl(h,h\bigr)_{s,\varphi}^{1/2}$.

The Hilbert space $H^{s,\varphi}(\Gamma)$ is separable,
continuously imbedded into the space $\mathcal{D}'(\Gamma)$, and
independent (up to equivalent norms) of the choice of the atlas
and the partition of unity. The collection of function spaces
$$
\{H^{s,\varphi}(\Gamma): s\in\mathbb{R},\;\varphi\in\mathcal{M}\}
\eqno(4)
$$
is naturally called the refined scale over the manifold $\Gamma$.

This scale admits an alternative intrinsic description. Let the
Riemannian structure on the manifold $\Gamma$ which defines the
density $dx$ be given (it is always possible), and let
$\triangle_{\Gamma}$ be the Beltrami-Laplace operator on $\Gamma$
. For $s\in\mathbb{R}$ and $\varphi\in\mathcal{M}$, we define the
function
$$
\varphi_{s}(t):=t^{s/2}\varphi(t^{1/2})\;\;\mbox{for}\;\;t\geq1\quad
\mbox{and}\quad \varphi_{s}(t):=\varphi(1)\;\;\mbox{for}\;\;0<t<1.
$$
We consider the operator $\varphi_{s}(1-\triangle_{\Gamma})$ in
the space $L_{2}(\Gamma,dx)$ as a Borel function of the
self-adjoint operator $1-\triangle_{\Gamma}$.

\begin{proposition}
For arbitrary $s\in\mathbb{R}$,
$\varphi\in\mathcal{M}$, the space $H^{s,\varphi}(\Gamma)$
coincides with the completion of the set of functions $u\in
C^{\infty}(\Gamma)$ with respect to the norm
$\|\varphi_{s}(1-\triangle_{\Gamma})\,u\|_{L_{2}(\Gamma,dx)}$
which is equivalent to the norm $\|u\|_{s,\varphi}$.
\end{proposition}

The following refinement of the classical Sobolev theorem
characterizes separating possibilities of scale (4).

\begin{proposition}
Let a function $\varphi\in\mathcal{M}$ and an
integer $\rho\geq0$ be given. The inequality
$$
\int_{1}^{\,+\infty}\frac{d\,t}{t\,\varphi^{\,2}(t)}<\infty
\eqno(5)
$$
is equivalent to the continuous imbedding
$H^{\rho+n/2,\varphi}(\Gamma)\hookrightarrow C^{\rho}(\Gamma)$.
The continuity of this imbedding implies its compactness.
\end{proposition}

\section{The basic results} We denote by $A^{+}$  a matrix
pseudodifferential operator formally adjoint to $A$ with respect
to the form $(\cdot,\cdot)_{\Gamma}$. We set
$$
N:=\left\{\,u\in\bigl(C^{\infty}(\Gamma)\bigr)^{p}:
\,Au=0\;\;\mbox{on}\;\;\Gamma\,\right\} \quad\mbox{and}\quad
N^{+}:=\left\{v\in\bigl(\,C^{\infty}(\Gamma)\bigr)^{p}:\,A^{+}v=0
\;\;\mbox{on}\;\;\Gamma\,\right\}.
$$
The ellipticity of system (2) implies that the spaces $N$ and
$N^{+}$ are finite-dimensional \cite[c. 52]{Agr90}.

\begin{theorem}
For each $s\in\mathbb{R}$,
$\varphi\in\mathcal{M}$ the linear bounded operator
$$
A:\,\prod_{k=1}^{p}\,H^{s+m_{k},\,\varphi}(\Gamma)\rightarrow
\bigl(H^{s,\varphi}(\Gamma)\bigr)^{p} \eqno(6)
$$
is defined.
It is a Fredholm one, has the kernel $N$ and the range
$$
\left\{f\in \bigl(H^{s,\varphi}(\Gamma)\bigr)^{p}:\,
\sum_{j=1}^{p}\,(f_{j},w_{j})_{\Gamma}=0\;\mbox{for each}\;
(w_{1},\ldots,w_{p})\in N^{+}\right\}.
$$
The index of the operator \rm (6) \it is equal to $\dim N-\dim
N^{+}$ and is independent of $s,\,\varphi$.
\end{theorem}

According to this theorem, $N^{+}$ is the defect subspace of
operator (6). Let's note \cite{AtSin63}, \cite[p. 32]{Agr90} that
in the scalar case ($p=1$) the index of operator (6) is equal to 0
if $\dim\Gamma\geq2$. Another sufficient condition for this
property is the ellipticity of system with a parameter on a
certain ray $K:=\{\lambda\in\mathbb{C}:\,\arg\lambda
=\mathrm{const}\} \, $\cite{Agr90}.

\begin{theorem}
For arbitrarily chosen parameters
$s\in\mathbb{R}$, $\varphi\in\mathcal{M}$ and $\sigma>0$, the
following a~priori estimate holds:
$$
\sum_{k=1}^{p}\;\bigl\|\,u_{k}\,\bigr\|_{s+m_{k},\varphi}\leq
c\,\left(\,\sum_{j=1}^{p}\;\bigl\|\,f_{j}\,\bigr\|_{s,\varphi}+
\sum_{k=1}^{p}\;\bigl\|\,u_{k}\,\bigr\|_{s-\sigma}\,\right).
$$
Here the number $c>0$ is independent of vector-functions $u$,
$f=Au$.
\end{theorem}

If the spaces $N$ and $N^{+}$ are trivial, then operator (6) is a
topological isomorphism. Generally, it is convenient to define the
isomorphism with the help of two projectors. Let's consider the
spaces in which operator (6) acts. Let us decompose them in the
following direct sums of the closed subspaces:
$$
\prod_{k=1}^{p}H^{s+m_{k},\,\varphi}(\Gamma)=N\dotplus\left\{u:\;
\sum_{k=1}^{p}\,(u_{k},v_{k})_{\Gamma}=0\;\mbox{for
each}\;(v_{1},\ldots,v_{p})\in N\right\},
$$
$$
\bigl(H^{s,\varphi}(\Gamma)\bigr)^{p}=N^{+}\dotplus
A\left(\,\bigl(H^{s,\varphi}(\Gamma)\bigr)^{p}\,\right).
$$
We denote by $P$ and $P^{+}$ respectively the oblique projectors
of these spaces onto the second terms in the sums. The projectors
are independent of $s,\varphi$.

\begin{theorem}
For arbitrary $s\in\mathbb{R}$, $\varphi\in\mathcal{M}$, the
restriction of operator $(6)$ onto the subspace
$P\left(\,\prod_{k=1}^{p}\,H^{s+m_{k},\,\varphi}(\Gamma)\,\right)$
establishes the topological isomorphism
$$
A:P\left(\,\prod_{k=1}^{p}\,H^{s+m_{k},\,\varphi}(\Gamma)\,\right)\leftrightarrow
P^{+}\left(\,\bigl(\,H^{s,\varphi}(\Gamma)\,\bigr)^{p}\,\right).
$$
\end{theorem}

\section{An application} Let $\Gamma_{0}$ be an open nonempty
subset of the manifold $\Gamma$. Denote
$$
H^{s,\varphi}_{\mathrm{loc}}(\Gamma_{0})
:=\left\{f\in\mathcal{D}'(\Gamma):\,\chi\,f\in
H^{s,\varphi}(\Gamma)\;\mbox{for each}\;\chi\in
C^{\infty}(\Gamma),\,\mathrm{supp}\,\chi\subset
\Gamma_{0}\right\}.
$$

\begin{theorem}
Suppose that a vector-function
$u\in\bigl(\mathcal{D}'(\Gamma)\bigr)^{p}$ is a solution of the
equation $Au=f$ on the set $\Gamma_{0}$, where $f\in
\bigl(\,H^{s,\varphi}_{\mathrm{loc}}(\Gamma_{0})\,\bigr)^{p}$ for
some parameters $s\in\mathbb{R}$ and $\varphi\in\mathcal{M}$. Then
$u\in\prod_{k=1}^{p}H^{s+m_{k},\varphi}_{\mathrm{loc}}(\Gamma_{0})$.
\end{theorem}

This theorem specifies, with regard to refined scale (4), known
propositions on local lifting of interior smoothness of an
elliptic system solution in the Sobolev scale (see e.g.
\cite{DN55, Her65, Ber65}). Note that the refined local smoothness
$\varphi$ of the right-hand side of the elliptic system is
inherited by its solution. Theorem 4 and Proposition 1 imply the
following sufficient condition for a chosen component $u_{k}$ of
the solution of system (2) to have continuous derivatives of a
prescribed order.

\begin{corollary}
Suppose that vector-functions
$u,f\in\bigl(\mathcal{D}'(\Gamma)\bigr)^{p}$ satisfy the equation
$Au=f$ on $\Gamma_{0}$. Let integers $\rho\geq0$, $k=1,\ldots,p$,
and a function $\varphi\in\mathcal{M}$ be such that inequality
$(5)$ is true. Then
$$
\Bigl(\;f_{j}\in
H^{\rho-m_{k}+n/2,\varphi}_{\mathrm{loc}}(\Gamma_{0})\;\;\mbox{for
every}\;j=1,\ldots,p\;\Bigr)\;\;\Rightarrow\;u_{k}\in
C^{\rho}(\Gamma_{0}).
$$
\end{corollary}

\end{document}